\documentclass{amsart}
\usepackage{amssymb}
\usepackage{latexsym}
\usepackage{amsmath}
\usepackage{euscript}
\usepackage{graphics}
\usepackage[all]{xy}

\newcommand{\be}{\begin{equation}}
\newcommand{\ee}{\end{equation}}
\newcommand{\ba}{\begin{eqnarray}}
\newcommand{\ea}{\end{eqnarray}}
\newcommand{\baa}{\begin{eqnarray*}}
\newcommand{\eaa}{\end{eqnarray*}}
\newcommand{\bb}{\begin{thebibliography}}
\newcommand{\eb}{\end{thebibliography}}

\newcommand{\bi}[1]{\bibitem{#1}}
\newcommand{\lab}[1]{\label{#1}}
\newcommand{\re}[1]{(\ref{#1})}

\newcounter{my}
\newcommand{\he}%
   {\stepcounter{equation}\setcounter{my}%
   {\value{equation}}\setcounter{equation}0%
   }%
\newcommand{\she}%
   {\setcounter{equation}{\value{my}}%
    }%

\renewcommand\t{\tilde}

\theoremstyle{definition}

\numberwithin{equation}{section}

\begin{document}

\title[Jordan algebras and orthogonal polynomials]
{Jordan algebras and orthogonal polynomials}

\author{Satoshi Tsujimoto}
\author{Luc Vinet}
\author{Alexei Zhedanov}

\address{Department of Applied Mathematics and Physics, Graduate School of Informatics,
Kyoto University, Sakyo-ku, Kyoto 606--8501, Japan}

\address{Centre de recherches math\'ematiques
Universite de Montr\'eal, P.O. Box 6128, Centre-ville Station,
Montr\'eal (Qu\'ebec), H3C 3J7}

\address{Institute for Physics and Engineering\\
R.Luxemburg str. 72 \\
83114 Donetsk, Ukraine \\}


\date{\today}

\begin{abstract}
We illustrate how Jordan algebras can provide a framework for the
interpretation of certain classes of orthogonal polynomials. The
big -1 Jacobi polynomials are eigenfunctions of a first order
operator of Dunkl type. We consider an algebra that has this
operator (up to constants) as one of its three generators and
whose defining relations are given in terms of anticommutators. It
is a special case of the Askey-Wilson algebra $AW(3)$. We show how
the structure and recurrence relations of the big -1 Jacobi
polynomials are obtained from the representations of this algebra.
We also present ladder operators for these polynomials and point
out  that the big -1 Jacobi polynomials satisfy the Hahn property
with respect to a generalized Dunkl operator.
\end{abstract}


\maketitle

\section{Introduction}
\setcounter{equation}{0} The relation between Lie groups and
special functions is a celebrated one and has been the object of
many books \cite{Vil}, \cite{Mil}, \cite{NSU}. Special functions
encode symmetries described by algebras. Over the years,
tremendous cross-fertilization has occurred within these areas.
The identification of new algebraic structures say like quantum
groups,  double affine Hecke algebras, polynomial algebras etc, in
connection with new manifestations of symmetry generally permit to
advance the theory of special functions. Conversely, the discovery
of new special functions has often revealed new algebraic tools.
This paper relates to the latter.

We have investigated recently the classes of orthogonal
polynomials that are eigenfunctions of Dunkl operators, that is of
differential-difference operators involving the reflection
operator.  These studies have allowed to discover new families of
classical polynomials.  They have also brought to light
interesting algebras defined through Jordan products, that is
anticommutators.

The big -1 Jacobi polynomials that we define next are the most
general eigenpolynomials of first order  Dunkl operators. They
have an associated Jordan algebra which is readily obtained  from
the operator of which they are the eigenfunctions. One purpose of
this paper is to illustrate  how the representations of this
algebra largely characterize the corresponding polynomials;
another is to provide the ladder operators for the big -1 Jacobi
polynomials.

Let us point out that the little -1 Jacobi polynomials have
already found applications in physical contexts. These polynomials
are obtained from the big -1 Jacobi polynomials by setting one
parameter equal to zero. They have been seen, in particular, to
arise in the wavefunctions of a supersymmetric Scarf Hamiltonian
with reflections \cite{PVZ1}.They have shown up furthermore in the
angular part of the separated solutions of the Schr\"odinger
equations associated to an infinite family of quantum
superintegrable models in the plane \cite{PVZ2}. This already
casts the Jordan algebra associated to the little -1 Jacobi
polynomials as the dynamical algebra of these physical systems.

\section{Big -1 Jacobi polynomials}
\setcounter{equation}{0}

The big -1 Jacobi polynomials $P_n^{(\alpha,\beta)}(x;c)$ are the
eigensolutions of  \be L^{(\alpha,\beta,c)}
P_n^{(\alpha,\beta)}(x;c) = \lambda_n P_n^{(\alpha,\beta)}(x;c),
\lab{LPP} \ee where the operator \be L^{(\alpha,\beta,c)}=
g_0(x)(R-I) + g_1(x) \partial_x R \lab{L_0} \ee with \be g_0(x)=
\frac{(\alpha+\beta+1)x^2 +(c\alpha -\beta)x + c}{x^{2}}, \quad
g_1(x)=\frac{2(x-1)(x+c)}{x}, \lab{g_01} \ee $I$ is the identity
operator and $R$ the reflection operator $Rf(x) = f(-x)$.

The eigenvalues of $L^{(\alpha,\beta,c)}$ are \be \lambda_n =
\left\{ {2n, \quad n \quad \mbox{even}    \atop
-2(\alpha+\beta+n+1), \quad n \quad \mbox{odd}} \right . .
\lab{lam-1} \ee As shown in \cite{VZ_D} the operator \re{L_0} is
the most general operator of the first order in   $\partial_x$ and
involving $R$ that has orthogonal polynomials as eigenfunctions.
The big -1 Jacobi polynomials depend on 3 parameters $\alpha,
\beta, c$. When $c=0$ they become the little -1 Jacobi polynomials
\cite{VZ_little}.

When $\alpha>-1, \: \beta>-1$ and $0<c<1$, the big -1 Jacobi
polynomials are orthogonal on the union of two symmetric intervals
$[-1,-c]$ and $[c,1]$ of the real axis: \be \int_{-1}^{-c}
P_n^{(\alpha,\beta)}(x;c) P_m^{(\alpha,\beta)}(x;c) w(x) dx +
\int_{c}^{1} P_n^{(\alpha,\beta)}(x;c) P_m^{(\alpha,\beta)}(x;c)
w(x) dx = h_n \: \delta_{nm} \lab{ort_P} \ee with positive
normalization constants $h_n$. The weight function is \be w(x) =
\frac{x}{|x|} (x+1)(x-c) (1-x^2)^{(\alpha-1)/2}(x^2
-c^2)^{(\beta-1)/2}. \lab{w_P} \ee

In \cite{VZ_little} and \cite{VZ_big}, we studied the properties
of the big and little -1 Jacobi polynomials directly from the
limit $q \to -1$ of the corresponding big and little q-Jacobi
polynomials.

We here wish to derive the structure relations and the recurrence
coefficients of these polynomials by starting from the eigenvalue
equation \re{LPP}. The main tool in this analysis will be a linear
Jordan algebra with 3 generators which is a special case of the
Askey-Wilson algebra $AW(3)$ \cite{Zhe} for $q=-1$. This algebra
was already presented in \cite{VZ_little}, \cite{VZ_big}; we
illustrate below its usefulness.

\section{Jordan algebra for big -1 Jacobi polynomials and intertwining operators}
\setcounter{equation}{0} Introduce the operators \be X=\frac{1}{2}
\left(L^{(\alpha,\beta,c)}+ \alpha+\beta+1  \right), \quad Y = 2x,
\quad Z = -\frac{2}{x}\left(c+(x-1)(x+c)R  \right). \lab{XYZ_def}
\ee In \cite{VZ_big} it was shown that these operators are closed
in frames of the Jordan algebra and that they satisfy the linear
anticommutation relations \be \{X,Y\} = Z+\omega_3, \quad
\{Y,Z\}=\omega_1, \quad \{Z,X\}=Y+\omega_2, \lab{anti_alg} \ee
where
$$
\omega_1=-8c, \; \omega_2 = 2(\alpha-\beta c), \; \omega_3
=2(\beta-\alpha c),
$$
with $\{A,B\} = AB + BA$ denoting as usual the anticommutator of
$A$ and $B$. (Note that our definition of the operators $X,Y$
slightly differs from that of \cite{VZ_big}.)

Observe that the first relation in \re{anti_alg} can be taken to
be the definition of $Z$ as the anticommutator (up to additive and
multiplicative constants) of the operator defining the eigenvalue
problem and the operator multiplication by $x$.

The Casimir operator of the algebra defined by \re{anti_alg} is
\be Q=Z^2+Y^2. \lab{Q_def} \ee In our realization the Casimir
operator takes the constant value $Q=4(c^2+1)$.

The algebra defined in \re{anti_alg} can be considered as the
limit $q  \to -1$ of the Askey-Wilson algebra $AW(3)$ \cite{Zhe}.
It also belongs to the class of Jordan algebras: the
anticommutators of any pair of generators are expressed in terms
of the generators. In the present case we have a linear Jordan
algebra with 3 generators $X,Y,Z$.

Consider  the canonical polynomial basis \be \Phi_n(x) = \left\{
{(x^2-c^2)^{\frac{n}{2}}, \quad n \quad \mbox{even}    \atop
(x+c)(x^2-c^2)^{\frac{n-1}{2}} , \quad n \quad \mbox{odd}} \right
. . \lab{basis_Phi} \ee The operator $L^{(\alpha,\beta,c)}$ is
2-diagonal and lower-triangular in this basis: \be
L^{(\alpha,\beta,c)} \Phi_n(x) = \lambda_n \Phi_n(x) + \eta_n
\Phi_{n-1}(x), \lab{LPhi} \ee where \be \eta_n= \left\{ {2(c-1)n,
\quad n \quad \mbox{even} \atop 2(c+1)(\beta+n), \quad n \quad
\mbox{odd}} \right . . \lab{eta_n} \ee Clearly, the operator $Y$
is 2-diagonal and upper-triangular in the basis $\Phi_n(x)$. The
existence of the basis $\Phi_n(x)$ with such properties is an
essential part of Terwilliger's approach to Leonard pairs
\cite{Ter}. From \re{LPhi} it is possible to obtain explicit
expressions for the big -1 Jacobi polynomials in terms of
hypergeometric functions \cite{VZ_big}.

An important property of the algebra \re{anti_alg} is that it
possesses simple intertwining operators $J_{\pm}$.  We define
these operators by the formulas \be
J_{+}=(Y+Z)(X-1/2)-\frac{\omega_2+\omega_3}{2} \lab{J+} \ee and
\be J_{-}=(Y-Z)(X+1/2)+\frac{\omega_2-\omega_3}{2}. \lab{J-} \ee
From the defining  relations \re{anti_alg}, we find that the
operators $J_{\pm}$ satisfy the anticommutation relations \be
\{X,J_{+}\} = J_{+}, \quad \{X,J_{-}\} = -J_{-}. \lab{com_JPM} \ee

It is seen readily that both $J_{+}^2$ and $J_{-}^2$ commute with
the operator $X$: \be [X,J_{+}^2]=[X,J_{-}^2]=0. \lab{com_X_JPM}
\ee   Let $\psi_n(x)$ be any eigenfunction of the operator $X$:
$$
X \psi_n(x) = \mu_n \psi_n(x),
$$
where \be \mu_n= \frac{1}{2}(\lambda_n + \alpha+\beta+1) =(-1)^n
\left(n+ \frac{\alpha+\beta+1}{2}  \right). \lab{mu} \ee From
\re{com_JPM} it is immediate to see  that the function $\t
\psi_n(x)=J_{+}\psi_n$ is also an eigenfunction of the operator
$X$ with the eigenvalue $\t \mu_n = 1-\mu_n$. In view of \re{mu}
$\t \mu_n = \mu_{n-1}$ if $n$ is even and $\t \mu_n = \mu_{n+1}$
if $n$ is odd.

Similar observations are made when $J_{-}$ replaces $J_{+}$.  The
function $J_{-} \psi_n$ is again an eigenfunction of the operator
$X$ with eigenvalue $-1-\mu_n$ which is $\mu_{n+1}$ for even $n$
and $\mu_{n-1}$ for odd $n$.

It is also clear that the operators $J_{\pm}$ transform
polynomials into polynomials (but not of the same degree).  This
means that the operator $J_{+}$ transforms big -1 Jacobi
polynomials of degree $n$ into big -1 Jacobi polynomials of degree
$n \mp 1$. Similarly, the operator $J_{-}$ transforms big -1
Jacobi polynomials of degree $n$ into big -1 Jacobi polynomials of
degree $n \pm 1$. (The upper sign corresponds to even $n$, the
lower sign to odd $n$.)

Taking into account the leading coefficients, we arrive at the
following formulas for the action of $J_{\pm}$ on big -1 Jacobi
polynomials  \be J_{+}P_n(x) = \left\{ { \frac{2(c-1)^2 n
(\alpha+\beta+n)}{\alpha+\beta+2n} P_{n-1}(x), \quad \mbox{if}
\quad n \quad \mbox{even}   \atop -2(\alpha+\beta+2(n+1))
P_{n+1}(x), \quad \mbox{if} \quad n \quad \mbox{odd}} \right .
\lab{J+_P} \ee and similarly \be J_{-}P_n(x) = \left\{
{2(\alpha+\beta+2(n+1)) P_{n+1}(x), \quad \mbox{if} \quad n \quad
\mbox{even}  \atop
-\frac{2(c+1)^2(\alpha+n)(\beta+n)}{\alpha+\beta+2n} P_{n-1}(x),
\quad \mbox{if} \quad n \quad \mbox{odd}} \right . . \lab{J-_P}
\ee Note that formulas \re{J+_P} and \re{J-_P} show that the
operators $J_{\pm}$ are block-diagonal in the basis of the
polynomials $P_n(x)$ with any block a $2\times 2$ matrix. They
basically provide a representation of the algebra \re{anti_alg}.

\section{Structure relations for the big -1 Jacobi polynomials}
\setcounter{equation}{0} We shall now derive the structure
relations for the big -1 Jacobi polynomials using the formulas
obtained in the previous section.

Define the operator $U_n^{(1)}, \: n=0,1,2,\dots$ as \be U_n^{(1)}
= \left\{ {J_{+}, \quad \mbox{if} \quad n \quad \mbox{even}  \atop
J_{-}, \quad \mbox{if} \quad n \quad \mbox{odd}} \right . .
\lab{U1} \ee
 We then have
 \be
 U_n^{(1)} P_n(x) =\epsilon^{(1)}_n P_{n-1}(x),  \lab{U1P} \ee
where \be \epsilon^{(1)}_n =\left\{{\frac{2(c-1)^2 n
(\alpha+\beta+n)}{\alpha+\beta+2n} , \quad \mbox{if} \quad n \quad
\mbox{even} \atop
-\frac{2(c+1)^2(\alpha+n)(\beta+n)}{\alpha+\beta+2n} , \quad
\mbox{if} \quad n \quad \mbox{odd}} \right . . \ee

Similarly, define the operator $U_n^{(2)}, \: n=0,1,2,\dots$ as
\be U_n^{(2)} = \left\{ {J_{-}, \quad \mbox{if} \quad n \quad
\mbox{even}  \atop J_{+}, \quad \mbox{if} \quad n \quad
\mbox{odd}} \right . . \lab{U2} \ee
 This entails
 \be
 U_n^{(2)} P_n(x) =\epsilon^{(2)}_n P_{n+1}(x), \lab{U2P} \ee
where \be \epsilon^{(2)}_n =\left\{{-2(\alpha+\beta+2(n+1)) ,
\quad \mbox{if} \quad n \quad         \mbox{even} \atop
2(\alpha+\beta+2(n+1)),  \quad \mbox{if} \quad n \quad \mbox{odd}}
\right . . \ee It is seen that the operator $U_n^{(1)}$ plays the
role of a lowering operator, while the operator  $U_n^{(2)}$ acts
as a raising operator for the polynomials $P_n(x)$. Relations such
as \re{U1P} and \re{U2P} are called structure relations: the
operators $U_n^{(1,2)}$ increment by $\pm 1$ the degree of the big
-1 Jacobi polynomials without altering the parameters $\alpha,
\beta, c$. Note however, that the operators $U_n^{(1,2)}$ depend
on $n$. Here this dependence is "minimal": even and odd sequences
of these operators do not depend on $n$.

These structure relations proved essential in demonstrating the
superintegrability of the two-dimensional infinite family of
quantum systems with Hamiltonians \cite{PVZ2}: \be
 H_k(r, \theta; \omega,
\alpha, \beta)=-\partial _r^2-\frac1r\partial _r-\frac
1{r^2}\partial_\theta^2 +\omega^2 r+\frac{\alpha
k^2}{2}\left(\frac{\frac\alpha 2-\cos k\theta\
R}{r^2\sin^2k\theta}\right)+\frac{\beta
k^2}{2}\left(\frac{\frac\beta 2-\sin k\theta\
R}{r^2\cos^2k\theta}\right), \lab{sup_H} \ee where $R$ is the
reflection operators with respect to $\theta$ \be
Rf(\theta)=f(-\theta). \lab{R_theta} \ee The following
restrictions on the real parameters $\alpha, \beta, $ and $k$:
$\alpha
>-1, \beta
> -1, $ $k \ne 0 $ are imposed with $r \in [0, \infty)$ and $-\frac{\pi}{2k} \leq \theta \leq
\frac{\pi}{2k}.$

The angular part of H (in the variable $\theta$) provides also a
supersymmetrization with reflections of the Scarf system in one
dimension \cite{PVZ1}. With the help of the above structure
relations, following the recurrence approach of \cite{KKM}, the
constants of motion of \re{sup_H} could be constructed
\cite{PVZ2}, thereby making the superintegrability manifest.
Furthermore, it was possible to determine the polynomial algebra
realized by these conserved quantities.

\section{Recurrence relations}
\setcounter{equation}{0}

Using the operators $J_{\pm}$ it is easy to derive the 3-term
recurrence relations of the  polynomials $P_n(x)$.

To that end we consider the operator \be
V=J_{+}(X+1/2)+J_{-}(X-1/2). \lab{V_def} \ee Recall that $X$ is
diagonal in the polynomials basis $P_n(x)$: $L P_n(x) = \mu_n
P_n(x)$, while $J_{\pm}$ are 2-diagonal in the same basis. Using
formulas \re{J+_P} and \re{J-_P} we have on the one hand, \be V
P_n(x) =\left\{{ \xi_n^{(-)} P_{n-1}(x) + \xi_n^{(+)} P_{n+1}(x),
\quad \mbox{if} \quad n \quad \mbox{even} \atop \eta_n^{(-)}
P_{n-1}(x) + \eta_n^{(+)} P_{n+1}(x) \quad \mbox{if} \quad n \quad
\mbox{odd}} \right . , \lab{V_2term} \ee where

$$
\xi_n^{(-)}= \frac{2(c-1)^2 \:
n(\alpha+\beta+n)(\mu_n+1/2)}{\alpha+\beta+2n}, \quad
\xi_n^{(+)}=2 (\mu_n-1/2) (\alpha+\beta+2(n+1))
$$
and
$$
\eta_n^{(-)}=-2(\mu_n-1/2) (\alpha+\beta+2(n+1)) , \quad
\eta_n^{(+)}=-
\frac{2(c+1)^2(\alpha+n)(\beta+n)(\mu_n+1/2)}{\alpha+\beta+2n}.
$$
(Recall that $\mu_n$ is the eigenvalue of the operator $X$ given
by \re{mu}.)

On the other hand, \re{J+} and  \re{J-} can be used to present
 the operator $V$ of \re{V_def} in the form \be V=2Y(X^2-1/4)
-\omega_3 X-\omega_2/2. \lab{V_def_1} \ee The operator $Y$
coincides with the operator multiplication by $2x$, so that in the
polynomial basis $P_n(x)$ we have $Y P_n(x) = 2x P_n(x)$. Hence
from \re{V_def_1} we have \be V P_n(x) =  (2(\mu_n^2-1/4)x
-\omega_3 \mu_n -\omega_2/2)P_n(x). \lab{VP_x} \ee Comparing
\re{V_2term} and \re{VP_x}, we arrive at the 3-term recurrence
relation for the polynomials $P_n(x)$: \be x P_n(x) = P_{n+1}(x) +
b_n P_n(x) + u_n P_{n-1}(x), \lab{3-term_P} \ee where the
recurrence coefficients are \be u_n = \left\{  {\frac{(1-c)^2
n(\alpha+\beta+n)}{(\alpha+\beta+2n)^2}, \quad n \quad \mbox{even}
\atop  \frac{(1+c)^2 (\alpha+n)(\beta+n)}{(\alpha+\beta+2n)^2},
\quad n \quad \mbox{odd} } \right . \lab{u_Jac} \ee and \be b_n =
\left\{  {-c +\frac{(c-1)n}{\alpha+\beta+2n}
+\frac{(1+c)(\beta+n+1)}{\alpha+\beta+2n+2}, \quad n \quad
\mbox{even} \atop  c + \frac{(1-c)(n+1)}{\alpha+\beta+2n+2} -
\frac{(c+1)(\beta+n)}{\alpha+\beta+2n}, \quad n \quad \mbox{odd} }
\right . . \lab{b_Jac} \ee The expressions \re{u_Jac}, \re{b_Jac}
coincide with the recurrence coefficients of the big -1 Jacobi
polynomials found in \cite{VZ_big} by a direct $q \to -1$ limit of
the big $q$-Jacobi polynomials counterparts.

\section{Lowering and raising operators for the big -1 Jacobi polynomials}
\setcounter{equation}{0} Consider the operator \be {\mathfrak D}
=A(x)(I-R) + B(x) \partial_x  + C(x)\partial_x R,  \lab{low_op}
\ee where \be A(x) = \frac{c^2}{x^3} -
\frac{c(c-1)}{2x^2}+\frac{\beta(c+1)}{2x}, \quad B(x) =
\frac{x^2-c^2}{x^2}, \quad C(x)=\frac{c(x+c)(1-x)}{x^2}. \lab{ABC}
\ee It is easily verified that the operator $\mathfrak D$
transforms any polynomial of degree $n$ into a polynomial of
degree $n-1$. Moreover, on the basis $\Phi_n(x)$ defined by
\re{basis_Phi} this operator acts simply as follows: \be
{\mathfrak D} \Phi_n(x)= \nu_n \: \Phi_{n-1}(x), \lab{D_Phi} \ee
where \be \nu_n = \left\{ {(1-c)n, \quad n \quad \mbox{even} \atop
(c+1)(\beta+n), \quad n \quad \mbox{odd} .} \right . . \lab{nu_n}
\ee When $c=0$ the operator  $\mathfrak D$ becomes the ordinary
Dunkl operator \be {\mathfrak D}\left |_{c=0} \right . =
\partial_x + \frac{\beta}{2x} (I-R). \lab{Dunkl_D} \ee Thus the
operator ${\mathfrak D}$ can be considered as a natural
generalization of the Dunkl operator with respect to the basis
$\Phi_n(x)$.

The operator ${\mathfrak D}$ satisfies an important intertwining property
\be
L^{(\alpha+2,\beta,c)} {\mathfrak D} + {\mathfrak D} L^{(\alpha,\beta,c)} + 2(\alpha+\beta+2) {\mathfrak D} =0. \lab{LD_inter} \ee
Relation \re{LD_inter} can be verified by direct calculations.

From \re{LD_inter} follows that for $\psi_n$ an eigenfunction of
the operator $L^{(\alpha,\beta,c)}$
$$
L^{(\alpha,\beta,c)} \psi_n(x) = \lambda_n \psi_n(x),
$$
the function $\tilde \psi_n(x) = {\mathfrak D} \psi_n(x)$ is an
eigenfunction of the operator $L^{(\alpha+2,\beta,c)}$: \be
L^{(\alpha+2,\beta,c)} \tilde \psi_n(x) = \tilde \lambda_n \tilde
\psi_n(x) \lab{Ltt} \ee with \be \t \lambda_n = -\lambda_n
-2(\alpha+\beta+2). \lab{tlam} \ee We note also that  the operator
${\mathfrak D}$ transforms polynomials of degree $n$ into
polynomials of degree $n-1$. Hence the operator ${\mathfrak D}$
transforms the big -1 Jacobi polynomials
$P_n^{(\alpha,\beta)}(x;c)$ (which are the unique polynomial
eigenfunctions of the operator $L^{(\alpha,\beta,c)}$) into the
polynomials $P_{n-1}^{(\alpha+2,\beta)}(x;c)$ (which are the
unique polynomial eigenfunctions of the operator
$L^{(\alpha+2,\beta,c)}$): \be {\mathfrak D}
P_n^{(\alpha,\beta)}(x;c) = \nu_n \:
P_{n-1}^{(\alpha+2,\beta)}(x;c), \lab{DPP} \ee where $\nu_n$ is
given by \re{nu_n}.

The operator ${\mathfrak D}$ is thus a lowering operator for the
big -1 Jacobi polynomials. Moreover, we proved that the big -1
Jacobi polynomials possess the Hahn property: namely that there
exists an operator ${\mathfrak D}$ such that its application to
the polynomials  $P_n^{(\alpha,\beta)}(x;c)$ gives another set of
orthogonal polynomials $P_n^{(\alpha+2,\beta)}(x;c)$.

In the case of little -1 Jacobi polynomials (i.e. when $c=0$), the
Hahn property had been proven in \cite{VZ_little}

It is well known \cite{Ger1} that the classical orthogonal
polynomials are completely characterized by the Hahn property with
respect to the ordinary derivative operator $\partial_x$. The
little -1 Jacobi polynomials satisfy their Hahn property with
respect to the Dunkl operator \re{Dunkl_D}. We see that the big -1
Jacobi polynomials satisfy the Hahn property with respect to a
generalized Dunkl operator ${\mathfrak D}$.  An interesting open
question is to characterize all orthogonal polynomials that have
the Hahn property with respect to generalized Dunkl operators
(which contain the operators $\partial_x$ and $R$).

Note also that the action of the operator ${\mathfrak D}$ is
equivalent to the application of two Christoffel transforms to the
big -1 Jacobi polynomials.

Indeed, it is well known (see, e.g. \cite{ZhS}) that the
Christoffel transform with parameter $a$ \be \t P_n(x) =
\frac{P_{n+1}(x) - \frac{P_{n+1}(a)}{P_n(a)} P_n(x)}{x-a}
\lab{CT_a} \ee is equivalent to the multiplication of their weight
function by a linear factor: $\t w(x) = (x-a) w(x)$. A multiple
Christoffel transform is hence equivalent to the multiplication of
the weight function by a polynomial: $\t w(x) = (x-a_1) (x-a_2)
\dots (x-a_N) w(x)$.

It is obvious from \re{w_P} that the weight function
$w^{(\alpha+2,\beta,c)}(x)$ corresponding to the big -1 Jacobi
polynomials  $P_n^{(\alpha+2,\beta)}(x;c)$ and the weight function
$w^{(\alpha,\beta,c)}(x)$ corresponding to the big -1 Jacobi
polynomials $P_n^{(\alpha,\beta)}(x;c)$, differ by a factor:
$1-x^2=(1-x)(1+x)$ \be w^{(\alpha+2,\beta,c)}(x) = (1-x^2)
w^{(\alpha,\beta,c)}(x). \lab{ww_CT} \ee Hence the big -1 Jacobi
polynomials ${\mathfrak D}P_n^{(\alpha,\beta)}(x;c)$ are obtained
from the polynomials $P_n^{(\alpha,\beta)}(x;c)$ by two successive
Christoffel transforms at the points $a_1=1$ and $a_2=-1$. It is
interesting to note that a similar property is valid for the
ordinary Jacobi polynomials $P_{n}^{(\alpha,\beta)}(x)$. Indeed,
Geronimus showed \cite{Ger1} that the polynomials $\partial_x
P_{n+1}^{(\alpha,\beta)}(x)$ are obtained from the Jacobi
polynomials $P_{n}^{(\alpha,\beta)}(x)$ by two successive
Christoffel transforms at the points $a_1=1$ and $a_2=-1$.

In the same manner one can construct a raising operator for the big -1 Jacobi polynomials.

Introduce the operator \be {\mathfrak R} = S_1(x)I + S_2(x) R +
T_1(x) \partial_x + T_2(x) \partial_x R , \lab{R_def} \ee where
\ba S_1(x)=&&\left( \beta\,c-\beta-2\,\alpha \right) x- \left( c+2
\right)
 \left( c-1 \right) + \nonumber \\&&{\frac {\beta-2\,{c}^{2}+2\,{c}^{2}\alpha-\beta\,
c}{x}}-{\frac {c \left( c-1 \right) }{{x}^{2}}}+{\frac
{2{c}^{2}}{{x }^{3}}}, \lab{S1} \ea \ba S_2(x)=&&\left(
\beta\,c-\beta+2\,c\alpha \right) x+ \left( c-1 \right)
 \left( 2\,c\alpha-c-2\,\beta \right) + \nonumber \\&&{\frac {-\beta-2\,{c}^{2}\alpha
+\beta\,c+2\,{c}^{2}}{x}}+{\frac {c \left( c-1 \right)
}{{x}^{2}}}- {\frac {{2c}^{2}}{{x}^{3}}} \lab{S2} \ea and \be
T_1(x) = {\frac { 2\left( x^2-1 \right)  \left( c^2-x^2
 \right)  }{{x}^{2}}}, \quad T_2(x) = {\frac {2c \left( 1+x \right)  \left( x-1 \right) ^{2} \left( x+c
 \right) }{{x}^{2}}}.
\lab{T1T2} \ee It is easily seen that the operator ${\mathfrak R}$
transforms any polynomial of degree $n$ into a polynomial of
degree $n+1$.

Moreover, there is an intertwining property \be
L^{(\alpha-2,\beta,c)} {\mathfrak R} + {\mathfrak R}
L^{(\alpha,\beta,c)} + 2(\alpha+\beta) {\mathfrak R} =0.
\lab{LR_inter} \ee Hence the operator $\mathfrak R$ maps the
polynomials $P_n^{(\alpha,\beta)}(x;c)$ to the polynomials
$P_{n+1}^{(\alpha-2,\beta)}(x;c)$: \be {\mathfrak R}
P_n^{(\alpha,\beta)}(x;c) = \kappa_n \:
P_{n+1}^{(\alpha-2,\beta)}(x;c), \lab{RPP} \ee where \be \kappa_n
= \left\{  { 2(c-1)(\alpha+\beta+n), \quad n \quad \mbox{even}
\atop  -2(c+1)(\alpha+n), \quad n \quad \mbox{odd} .} \right . .
\lab{kappa_n} \ee When $c=0$, the expression for the raising
operator becomes more simple: \be {\mathfrak R}=2(1-x^2)
\partial_x R - \frac{\beta(x-1)^2}{x} R +(2 +\frac{\beta}{x}-\left(\beta+2
\alpha)x\right)I. \lab{R_l} \ee Expression \re{R_l} gives the
raising operator for the little -1 Jacobi polynomials found in
\cite{VZ_little}.

We mentioned in the introduction that the little -1 Jacobi
polynomials form the wavefunctions of the supersymmetric Scarf
Hamiltonian (which is the angular part of \re{sup_H}). As was
pointed out in \cite{PVZ1}, in view of their intertwining
relations, for $c=0$, the operators \re{Dunkl_D} and \re{R_l}
allow to connect eigenfunctions of such extended Scarf
Hamiltonians with different parameters.

\section{Conclusions}

We have shown that algebras, whose defining relations are given in
terms of Jordan products, can be naturally associated to certain
families of orthogonal polynomials. We have used for purposes of
illustration the example of the big and little -1 polynomials. We
made the case of the usefulness of such Jordan algebras in
deriving structural properties of the corresponding polynomials.
We noted that there are some physical systems for which the
algebra we have studied is dynamical. We developed in \cite{BI}
the theory of the Bannai-Ito polynomials where again a Jordan
algebra was key. We trust in concluding that such structures will
continue to appear in various guises (physical and mathematical)
and will warrant further analysis.

\bigskip\bigskip
{\Large\bf Acknowledgments}
\bigskip

\noindent The authors would like to gratefully acknowledge the
hospitality extended to LV and AZ by Kyoto University and to ST
and LV by the Donetsk Institute for Physics and Technology  in the
course of this investigation. The research of ST is supported in
part through funds provided by  KAKENHI (22540224), JSPS. The
research of LV is supported in part by a research grant from the
Natural Sciences and Engineering Research Council (NSERC) of
Canada.

\newpage

\bb{99}

\bi{AA} G.E.Andrews and R.Askey, {\it Classical orthogonal
polynomials}, Polyn\^omes  Orthogonaux et Applications, Lecture
Notes in Mathematics, 1985, V. {\bf 1171}, 36--62.






\bi{BI} E. Bannai and T. Ito, {\it Algebraic Combinatorics I:
Association Schemes}. 1984. Benjamin \& Cummings, Mento Park, CA.

\bi{Cheikh} Y. Ben Cheikh and M.Gaied, {\it Characterization of the Dunkl-classical symmetric orthogonal polynomials}, Appl. Math. and Comput.
{\bf 187}, (2007) 105--114.


\bi{Chi} T. Chihara, {\it An Introduction to Orthogonal
Polynomials}, Gordon and Breach, NY, 1978.

\bi{Chi_Chi} L.M.Chihara,T.S.Chihara, {\it A class of nonsymmetric orthogonal polynomials}. J. Math. Anal. Appl. {\bf 126} (1987), 275--291.

\bi{Dunkl} C.F.Dunkl, {\it Integral kernels with reflection group invariance}. Canadian Journal of Mathematics, {\bf 43} (1991)
1213--1227.


\bi{Ger1} Ya.L.Geronimus, {\it On polynomials orthogonal with respect to to
the given numerical sequence and on Hahn's theorem}, Izv.Akad.Nauk, {\bf 4}
(1940), 215--228 (in Russian).

\bi{GVA} J.Geronimo and W. Van Assche, {\it Orthogonal polynomials on several intervals via a polynomial mapping},
Trans. Amer. Math. Soc., {\bf 308}(2) (1988), 559--581.








\bi{IT} T.Ito, P.Terwilliger, {\it Double Affine Hecke Algebras of Rank 1
and the $Z_3$-Symmetric Askey-Wilson Relations}, SIGMA {\bf 6} (2010), 065, 9 pages.

\bi{KS} R.Koekoek, R.Swarttouw, {\it The Askey-scheme of hypergeometric orthogonal polynomials and its q-analogue}, Report no. 98-17, Delft University of Technology, 1998.

\bi{KLS} R. Koekoek,P. Lesky, R. Swarttouw, {\it Hypergeometric Orthogonal Polynomials and Their Q-analogues}, Springer-Verlag, 2010.

\bi{K_Racah} T.Koornwinder, {\it On the limit from q-Racah polynomials to big q-Jacobi polynomials}, ArXiv:1011.5585.

\bi{KKM} E.K.Kalnins, J.M.Kress, W.Miller, {\it A Recurrence
Approach to higher order quantum superintegrability}, SIGMA {\bf
7} (2011), 031.


\bi{Kui} A.B.J. Kuijlaars, A. Martinez-Finkelshtein, and R. Orive, {Orthogonality of Jacobi polynomials with general
parameters}, Electronic Trans. Numer. Anal. {\bf 19} (2005), 1--17.

\bi{Leonard} D.Leonard, {\it Orthogonal Polynomials, Duality and
Association Schemes}, SIAM J. Math. Anal. {\bf 13} (1982)
656--663.


\bi{MP} F.Marcell\'an and J.Petronilho, {\it Eigenproblems for Tridiagonal 2-Toeplitz Matrices
and Quadratic Polynomial Mappings},  Lin. Alg. Appl. {\bf 260} (1997)   169--208.

\bi{Mil} W.Miller, {\it Lie Theory and Special Functions},
Academic Press, 1968.



\bi{NSU} A.F. Nikiforov, S.K. Suslov, and V.B. Uvarov, {\em
Classical Orthogonal Polynomials of a Discrete Variable},
Springer, Berlin, 1991.


\bi{PVZ1} S.Post, L.Vinet and A.Zhedanov, {\it Supersymmetric
Quantum Mechanics with Reflections}, arXiv:1107.5844.

\bi{PVZ2} S.Post, L.Vinet and A.Zhedanov, {\it An infinite family
of superintegrable Hamiltonians with reflection in the plane}, in
preparation.


\bi{Sz} G. Szeg\H{o}, Orthogonal Polynomials, fourth edition,  AMS, 1975.

\bi{Ter} P. Terwilliger, {\it Two linear transformations each tridiagonal with respect to an eigenbasis of the other}.
Linear Algebra Appl. {\bf 330} (2001) 149--203.

\bi{Ter2} P.Terwilliger, {\it Two linear transformations each
tridiagonal with respect to an eigenbasis of the other; an
algebraic approach to the Askey scheme of orthogonal polynomials},
arXiv:math/0408390v3.

\bi{BI} S.Tsujimoto, L.Vinet and A.Zhedanov, {\it Dunkl shift
operators and Bannai-Ito polynomials}, arXiv:1106.3512.

\bi{Vil} N. Ja. Vilenkin  and A. U. Klimyk, {Representation of Lie
groups and special functions}, Volumes 1--3, Math. Appl., Kluwer
Acad. Publ., Dordrecht, 1991-1992.



\bi{VZ_little} L.Vinet and A.Zhedanov, {\it A ``missing`` family
of classical orthogonal polynomials}, arXiv:1011.1669v2.

\bi{VZ_D} L.Vinet and A.Zhedanov, {\it A Bochner Theorem for Dunkl
Polynomials}, SIGMA {\bf 7} (2011), 020, 9 pages.

\bi{VZ_big} L.Vinet and A.Zhedanov, {\it A limit $q=-1$ for big q-Jacobi polynomials},
Trans.Amer.Math.Soc., to appear,    arXiv:1011.1429v3.


\bi{Zhe} A. S. Zhedanov. {\it Hidden symmetry of Askey-Wilson polynomials}, Teoret. Mat. Fiz.
{\bf 89} (1991) 190--204. (English transl.: Theoret. and Math. Phys. {\bf 89} (1991), 1146--1157).

\bibitem{ZhS} A.S. Zhedanov, {\it Rational spectral transformations
and orthogonal polynomials}, J. Comput. Appl. Math. {\bf 85}, no. 1
(1997), 67--86.

\eb

\end{document}